\documentclass[letterpaper,journal]{IEEEtran}
\usepackage{amsmath,amsfonts}
\usepackage{algorithmic}
\usepackage{algorithm}
\usepackage{array}
\usepackage[caption=false,font=normalsize,labelfont=sf,textfont=sf]{subfig}
\usepackage{textcomp}
\usepackage{stfloats}
\usepackage{url}
\usepackage{verbatim}
\usepackage{graphicx}
\usepackage{subcaption}
\usepackage{cite}
\usepackage{amssymb}
\usepackage{booktabs}
\usepackage{multirow}
\usepackage{xcolor}
\usepackage{cleveref}
\usepackage{optidef}
\usepackage{flushend}
\newcommand{\R}{\mathbb{R}}
\newcommand{\E}{\mathbb{E}}
\newcommand{\Prob}{\mathbb{P}}

\begin{document}

\title{Reinforcement Learning for Public Safety Power Shutoffs Under\\
Decision-Dependent Uncertainty and Nonlinear Wildfire Ignition Models}

\author{Prasanna Raut,~\IEEEmembership{Student~Member,~IEEE,} Chaoyue Zhao, \IEEEmembership{Member,~IEEE,} Alexandre Moreira, \IEEEmembership{Senior Member,~IEEE} }

\maketitle

\begin{abstract}
Power grid infrastructure is an increasingly significant source of wildfire ignitions and poses severe risks to communities in fire-prone regions. Public Safety Power Shutoffs (PSPS) have emerged as a critical operational tool for utilities to mitigate this risk by proactively de-energizing portions of the grid under high-threat conditions. These shutoffs, however, impose costs on affected communities, and it is therefore essential that PSPS decisions be informed by realistic models of wildfire ignition risk. Current Mixed Integer Programming based methods require restrictive structural assumptions about the probability models for line failures caused by power line ignitions. While these simplifications yield tractable solutions, the resulting models may differ significantly from the true underlying dynamics. In this paper, we propose a reinforcement learning framework based on Proximal Policy Optimization that learns to adjust the topology of a distribution system by interacting directly with a simulator that accommodates any line failure probability model without imposing such restrictions. We test our methodology on 54-bus and 138-bus distribution systems and demonstrate its ability to lower operational costs compared to existing methods while allowing only marginally increased compute times as network size grows.
\end{abstract}


\section{Introduction}
\label{sec:introduction}
 
Wildfires pose an existential threat to communities in fire-prone regions, and power grid infrastructure has been increasingly identified as a significant ignition source~\cite{cpuc2026wildfire}. High-profile events underscore the severity of this risk. In 2023, downed power lines sparked fires in Lahaina, Hawaii\cite{mauicounty2024lahaina}. In February 2024, a decayed utility pole ignited the largest wildfire in Texas history on the panhandle \cite{texaswildfirereport2024}. The 2018 Camp Fire, which is the deadliest and most destructive wildfire in California history, was caused by a single faulty hook on a PG\&E transmission line~\cite{butteda2020campfire}. As climate change intensifies drought and heat, the frequency and severity of such interactions are expected to grow~\cite{cpuc2026wildfire}. Beyond the immediate fire risk, failures in grid infrastructure can propagate into large-scale cascading outages. Research on the North American power grid shows that even a small, topologically central set of components is disproportionately responsible for triggering large cascades~\cite{yang2017small}, which indicates that a wildfire-induced equipment failure can threaten the stability of the broader network.

Public Safety Power Shutoffs (PSPS) have emerged as the primary operational tool for utilities to manage this risk by proactively de-energizing portions of the grid under high-threat conditions~\cite{cpuc2026psps}. Following their widespread deployment by California utilities since 2019, PSPS events are now being adopted by utilities across the country. For example, Idaho Power and Xcel Energy in Colorado each implemented shutoffs for the first time in 2024~\cite{cpuc2026psps}. However, PSPS events impose significant costs on affected communities. De-energizing lines reduces wildfire ignition risk but simultaneously causes load loss, with documented disparate impacts on medically vulnerable individuals, environmental justice communities, and low-income households~\cite{sotolongo2020california}. Designing PSPS strategies that effectively balance wildfire risk against the costs of outages is therefore a critical operational problem.

The research community has addressed this balance through a range of optimization models. The work in~\cite{trakas2017optimal} proposes a stochastic programming framework for improving distribution system resilience against approaching wildfires. By accounting for dynamic line ratings and uncertainties in solar radiation, wind speed, and wind direction, their method is able to attain higher state-of-charge (SOC) in the batteries, which then serves the demand for the system when the connection to the up-stream grid is lost due to the wildfire. At the operational level,~\cite{yang2022resilient} introduces the use of microgrids to manage wildfire risks without resorting to outages, validating the approach on a practical grid scenario from the 2019--2020 Australian wildfire season. The optimal power shutoff problem is formulated in~\cite{rhodes2020balancing}, where the objective is to maximize delivered power while minimizing ignition risk by selectively de-energizing grid components. This line of work is extended in~\cite{su2023quasi} through a quasi second-order stochastic dominance model that incorporates mobile power sources to reduce the load shedding imposed on customers during PSPS events. On the investment planning side,~\cite{pianco2025decision} proposes a distributionally robust methodology to identify optimal combinations of new lines, hardened infrastructure, and switching devices that minimize the long-term operating cost of the distribution system, thereby mitigating the risk of grid-caused wildfires through targeted infrastructure investment.

A distinguishing feature shared by the more recent operational models is the explicit treatment of decision-dependent uncertainty (DDU), i.e. the observation that the probability of a line failure is not determined solely by exogenous factors such as weather and vegetation, but also depends on the operator's own dispatch and switching decisions. When power flows through a wildfire-prone line, the thermal stress on conductors increases the probability of a fault or ignition event~\cite{yang2025tree}. The work in~\cite{moreira2024distribution} formalizes this structure within a two-stage distributionally robust optimization framework, in which switching actions in the first stage shape the ambiguity sets governing line failure probabilities. The same operational setting is addressed in~\cite{zhang2025power} using a chance-constrained stochastic programming model that incorporates spatial correlations across line failures via an Archimedean copula and models the decision-dependent impact of power flows through a distortion function applied to the failure probability distributions.

Despite their contributions, both~\cite{moreira2024distribution} and~\cite{zhang2025power} rely on linear or piecewise-linear models for the relationship between power flow and line failure probability. Empirical evidence suggests, however, that the true relationship may exhibit threshold effects. The work in~\cite{Muhs2021linefailure} characterizes the probability of wildfire ignition caused by distribution lines and finds that ignition risk can experience sharp jumps at certain power flow levels, and~\cite{van2019fire} reports massive reductions in ignition risk when flows are constrained below critical thresholds. Because the true failure probability model is difficult to capture through physics-based methods, it is better learned from data, and the resulting curve may be linear, piecewise-linear, step-shaped, or anything in between. We therefore seek a methodology that is flexible enough to accommodate any such curve. A further limitation is computational in nature. Both methods require solving mixed-integer programs (MIPs), which have an exponential worst-case time complexity in terms of the network size. Furthermore, handling DDU in the optimization framework is challenging and can lead to a non-linear problem reformulation and linearization of it requires introduction of additional integer variables and constraints, which further increases the size of the MIPs.

To address these limitations, we formulate the PSPS problem as a Markov Decision Process (MDP). The network's operational configuration is represented as a state defined by the current topology and system status, while actions correspond to switching decisions that energize or de-energize selected lines. State transitions capture both exogenous wildfire-related conditions and the endogenous influence of power flows on line failure probabilities. The objective is to minimize total operational cost over a wildfire event, including energy procurement costs, unserved energy penalties, and switching costs. Our approach is based on Proximal Policy Optimization (PPO)~\cite{schulman2017proximalpolicyoptimizationalgorithms}, a deep reinforcement learning algorithm, to solve this MDP. This approach offers two key advantages over optimization-based methods. First, RL algorithms learn policies through interaction with a simulator, requiring only the ability to sample outcomes from the failure probability distribution rather than an analytical representation. This enables incorporation of arbitrarily complex failure models, including step threshold functions and data-driven models learned from historical outage records without any structural constraints imposed by tractability requirements. Second, by using RL to determine switch positions, we avoid solving computationally expensive MIPs.

Deep reinforcement learning has been applied to a range of power grid resilience problems~\cite{cao2025deep}. At the distribution level,~\cite{dehghani2021intelligent} develops a deep RL planning framework that selects hardening strategies to enhance long-term resilience against hurricane-caused component failures, demonstrating improvements of over 30\% relative to standard engineering strategies. The work in~\cite{hosseini2021resilient} develops a resilience controller that learns to dispatch distributed generation and energy storage during a progressing hurricane, modeling distribution grid operation under uncertainty as an MDP. For wildfire settings specifically,~\cite{kadir2024RL} proposes a deep RL decision-support system for power grid operation that provides generation setpoints and line switching actions to minimize outages for a transmission grid during a progressing wildfire. A common feature of these methods is that the wildfire or hurricane is treated as an exogenous source of uncertainty, in the sense that the grid operator responds to an externally evolving hazard but does not influence the hazard itself. Our setting is fundamentally different. Power flows through wildfire-area lines affect the probability of ignition, making the uncertainty endogenous to the operator's decisions. Beyond this decision-dependent structure, the PSPS setting poses additional challenges for deep RL. Depending on the network size and number of switchable lines, the set of admissible topologies can become enormous, and deep RL is known to struggle in problems with large action spaces~\cite{dulac2015deep}, while the heavy-tailed magnitude of operating costs (and hence rewards) can further impede learning~\cite{vanHasselt2016rewardnormalization}. Addressing these issues jointly requires modeling and algorithmic choices that none of the above frameworks provide.

We summarize the main contributions of this paper below:
\begin{enumerate}
    \item A PPO-based reinforcement learning framework for PSPS decision-making that accommodates nonlinear line failure probability models through simulator-based learning, with a methodology for continuous action space discretization that guarantees radiality-feasible topology selections. The framework further incorporates per-trajectory reward standardization to stabilize training under the heavy-tailed, heterogeneous cost distributions that arise in the PSPS setting.
    \item Numerical evaluation on 54-bus and 138-bus distribution systems demonstrating two complementary consequences of linear model mismatch: risk overestimation at sub-threshold flows, which causes optimization-based baselines to perform unnecessary switching that the PPO agent correctly avoids; and risk underestimation under an extreme wildfire scenario, where the baselines leave dangerous flows on wildfire-area lines that the PPO agent eliminates through load isolation, achieving substantially fewer failures and lower costs. 
\end{enumerate}

The remainder of this paper is organized as follows. Section~\ref{sec:problem} describes the distribution system model, MDP formulation, and operational constraints. Section~\ref{sec:methodology} presents the PPO-based solution algorithm. Section~\ref{sec:case_study} examines the effect of failure model nonlinearity on policy performance. Section~\ref{sec:case_study_138} evaluates policy behavior under an extreme wildfire scenario. Section~\ref{sec:conclusion} concludes with directions for future work.

\section{Problem Description}
\label{sec:problem}

In this paper, we develop an MDP approach to address distribution system operation under wildfire risk considering decision-dependent uncertainty. We describe our problem formulation in the next subsections.

\subsection{Distribution System Model Setup}

Consider a distribution system represented by a graph $\mathcal{G} = (\mathcal{N}, \mathcal{L})$, where $\mathcal{N}$ denotes the set of buses and $\mathcal{L}$ the set of lines. The bus set is partitioned into substation buses $\mathcal{N}^{\text{sub}}$ that inject power from the transmission system, and load buses $\mathcal{N} \setminus \mathcal{N}^{\text{sub}}$ with active and reactive demand requirements.

A subset $\mathcal{L}^{\text{sw}} \subseteq \mathcal{L}$ are switchable lines whose energization status can be controlled by the operator. A subset $\mathcal{L}^{\text{fr}} \subseteq \mathcal{L}$ traverses wildfire-prone areas where energized lines may fail with a higher than normal probability that depends on operating conditions. We use $av_l \in \{0,1\}$ to denote the availability of line $l$, where $av_l = 0$ indicates failure that cannot be restored within the horizon.
\subsection{MDP Formulation}

We formulate the PSPS problem as a MDP defined by the tuple $(\mathcal{S}, \mathcal{A}, P, R)$.

\subsubsection{State Space}

The state at time $t$ is:
\begin{equation}
    \boldsymbol{s}_t = \left[\boldsymbol{av}_t^\top,\; \boldsymbol{z}_t^{\text{sw},0\top},\; \boldsymbol{D}_t^{p\top},\; \boldsymbol{D}_t^{q\top}\right]^\top
\end{equation}
where $\boldsymbol{av}_t = [av_{l,t}]_{l \in \mathcal{L}} \in \{0,1\}^{|\mathcal{L}|}$ is the line availability vector, $\boldsymbol{z}_t^{\text{sw},0} = [z_{l,t}^{\text{sw},0}]_{l \in \mathcal{L}^{\text{sw}}} \in \{0,1\}^{|\mathcal{L}^{\text{sw}}|}$ is the pre-decision switch status (1 = closed), and $\boldsymbol{D}_t^p, \boldsymbol{D}_t^q \in \R^{|\mathcal{N}|}$ are the active and reactive demand vectors at all buses.

The line availability $\boldsymbol{av}_t$ and switch status $\boldsymbol{z}_t^{\text{sw},0}$ are endogenous states that evolve based on operator decisions and stochastic failures. The demand profiles $\boldsymbol{D}_t^p, \boldsymbol{D}_t^q$ are exogenous parameters that evolve independently.

\subsubsection{Actions}

At each stage, the primary control action for the MDP is the choice of switching configuration:
\begin{equation}
    \boldsymbol{z}_t^{\text{sw}} = \left[z_{l,t}^{\text{sw}}\right]_{l \in \mathcal{L}^{\text{sw}}} \in \{0,1\}^{|\mathcal{L}^{\text{sw}}|} 
    \label{eq:action}
\end{equation}
where $z_{l,t}^{\text{sw}} = 1$ indicates that line $l$ is closed (energized) in the new configuration. The switching operation indicator is $y_{l,t}^{\text{sw}} = |z_{l,t}^{\text{sw}} - z_{l,t}^{\text{sw},0}|$, which equals 1 if the status of line $l$ changes at time $t$.

Once switch positions are fixed, the operator also solves for continuous operating variables including active and reactive power flows $f_l^p, f_l^q$ on each line, squared bus voltages $v_b$, active and reactive substation injections $p_b^{\text{sub}}, q_b^{\text{sub}}$, the active and reactive load surplus $\Delta D_b^{{p}+}, \Delta D_b^{{q}+}$ and the active and reactive load shedding variables $\Delta D_b^{{p}-}, \Delta D_b^{{q}-}$. These continuous variables are determined by solving the operational optimization problem described in Section~\ref{sec:opf}, and thus constitute secondary actions that are not part of the MDP policy.

\subsubsection{Transition Dynamics and Decision-Dependent Uncertainty}

State transitions are driven by the stochastic evolution of line availability $av_{l,t}$. We model the next-period availability of each available line $l \in \mathcal{L}$ as a Bernoulli random variable whose failure probability depends on power flow:
\begin{equation}
    \Prob\!\left(av_{l,t+1} = 0 \;\middle|\; av_{l,t} = 1,\; |f_{l,t}^p|\right) =
        g_l\!\left(|f_{l,t}^p|\right)  
    \label{eq:transition}
\end{equation}
Once a line has failed, it remains unavailable for the remainder of the horizon: $\Prob(av_{l,t+1} = 1 \mid av_{l,t} = 0) = 0$.

Because $f_{l,t}^p$ belongs to the stage-$t$ action, \eqref{eq:transition} couples switching and dispatch decisions to future network availability; this is the source of decision-dependent uncertainty in the problem. The operator's choice of topology affects not only the current cost, but also the distribution of future line failures.

\subsubsection{Nonlinear Failure Probability Model}

The function $g_l(\cdot)$ in \eqref{eq:transition} specifies how the failure probability of a line depends on the power flow it carries. The linear model used in prior work~\cite{moreira2024distribution} is:
\begin{equation}
    g_l^{\text{linear}}\!\left(|f_l^p|\right) = \gamma_l + \beta_l \cdot |f_l^p|
    \label{eq:linear_model}
\end{equation}
where $\gamma_l \in (0,1]$ is the baseline failure probability and $\beta_l \geq 0$ captures the sensitivity to power flow.

Empirical evidence suggests~\cite{Muhs2021linefailure}, however, that failure risk may be low for lightly loaded lines but increase sharply once flows exceed a critical thermal threshold. To model this behavior, we consider a \textbf{step threshold model}:
\begin{equation}
    g_l^{\text{step}}\!\left(|f_l^p|\right) =
    \begin{cases}
        \gamma_l + \beta_l \cdot F_l^{\max} & \text{if } |f_l^p| > \tau \cdot F_l^{\max} \\
        \gamma_l & \text{otherwise}
    \end{cases}
    \label{eq:step_model}
\end{equation}
where $\tau \in [0,1]$ is the threshold fraction and $F_l^{\max}$ is the line's thermal capacity. When $|f_l^p|$ is below $\tau \cdot F_l^{\max}$, the line operates at its baseline failure probability; once the flow exceeds the threshold, the failure probability increases by $\beta_l \cdot F_l^{\max}$.

\begin{figure}[t]
    \centering
    \includegraphics[width=\linewidth]{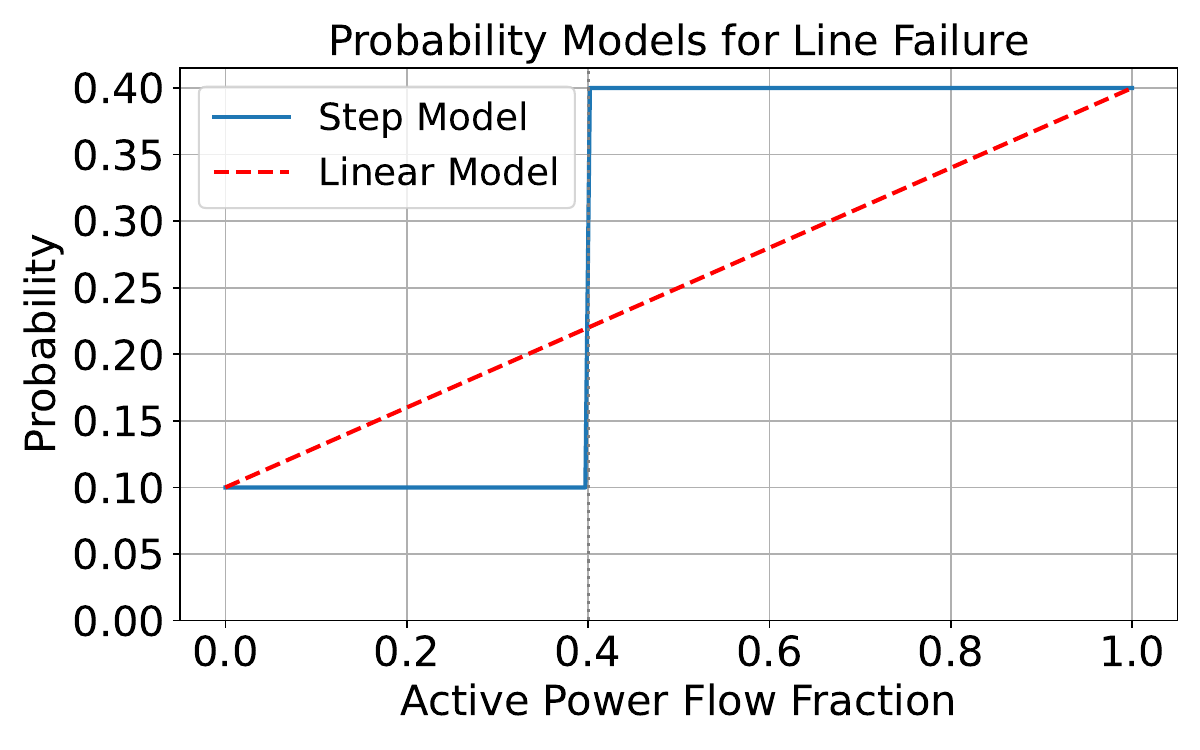}
    \caption{Comparison of line failure probability models as a function of power flow (normalized to line capacity). The static optimization-based DDU method assumes a linear relationship (dashed), while the true model may exhibit step threshold behavior (solid) where failure risk increases sharply only above a critical loading level.}
    \label{fig:ignition_models}
\end{figure}

\Cref{fig:ignition_models} illustrates the difference between the linear and step threshold models. The key implication is that the linear model continuously penalizes any positive flow, whereas the step model implies that moderate flows carry negligible additional risk---a behavior that may lead to qualitatively different optimal policies.

Note: While we have used the step threshold model for clarity, our method applies for any arbitrary failure probability function $g_l(\cdot)$ in \eqref{eq:transition}.

\subsubsection{Reward Function}

The immediate reward at time $t$ is the negative operating cost:
\begin{align}
    r_t = &-\sum_{b \in \mathcal{N}^{\text{sub}}} C^{\text{energy}} \cdot p_{b,t}^{\text{sub}}
           - \sum_{l \in \mathcal{L}^{\text{sw}}} C^{\text{switch}} \cdot y_{l,t}^{\text{sw}} \notag \\
          &- \sum_{b \in \mathcal{N}} C^{\text{load\_loss}} \cdot
             \Bigl(\Delta D_{b,t}^{p+} + \Delta D_{b,t}^{p-} \notag \\
          & \hspace{65pt} + \Delta D_{b,t}^{q+} + \Delta D_{b,t}^{q-}\Bigr)
    \label{eq:reward}
\end{align}
where $C^{\text{energy}}$ is the cost per unit of active power procured from substations, $C^{\text{switch}}$ is the per-operation switching cost, and $C^{\text{load\_loss}}$ is the penalty per unit of load not served (or over-supplied).

\subsubsection{State-Action Value Function}

The objective is to find a policy $\pi: \mathcal{S} \to \Delta(\mathcal{A})$ that maximizes the expected discounted cumulative reward:
\begin{equation}
    J_T(\pi) = \E_\pi\!\left[\sum_{t=1}^{T} r_t\right]
\end{equation}
where $T$ is the number of decision epochs. The state-value function and state-action value function are:
\begin{align}
    V_{T}^\pi(\boldsymbol{s}) &= \E_\pi\!\left[\sum_{t=1}^{T} r_{t} \;\middle|\; \boldsymbol{s}_0 = \boldsymbol{s}\right] \\
    Q_{T}^\pi(\boldsymbol{s}, \boldsymbol{a}) &= \E_\pi\!\left[\sum_{t=1}^{T} r_{t} \;\middle|\; \boldsymbol{s}_0 = \boldsymbol{s},\, \boldsymbol{a}_1 = \boldsymbol{a}\right]
\end{align}

\subsection{Operational Constraints}
\label{sec:opf}

Once the switching actions $\boldsymbol{z}_t^{\text{sw}}$ \eqref{eq:action} of the MDP are decided, the topology of the distribution system is fixed. The resulting power flows, power injections, unserved and surplus load must satisfy the operational constraints of the physical infrastructure. These quantities feed into the transition dynamics \eqref{eq:transition} and the reward \eqref{eq:reward}. The operator determines them by solving the following LP, which uses the linearized AC formulation of~\cite{mashayekh2017linearACpf, moreira2024distribution}.
\begin{align}
    \min_{\substack{p_{b,t}^{\text{sub}},\, q_{b,t}^{\text{sub}},\, f_{l,t}^{p},\, f_{l,t}^{q}, \\[2pt] v_{b,t},\, \Delta D_{b,t}^{p\pm},\, \Delta D_{b,t}^{q\pm}}} & \quad
     \sum_{b \in \mathcal{N}^{\text{sub}}} C^{\text{energy}} \cdot p_{b,t}^{\text{sub}}
    + \sum_{l \in \mathcal{L}^{\text{sw}}} C^{\text{switch}} \cdot y_{l,t}^{\text{sw}} \notag \\
     + \sum_{b \in \mathcal{N}} C^{\text{load\_loss}} & \cdot
    \Bigl(\Delta D_{b,t}^{p+} + \Delta D_{b,t}^{p-} + \Delta D_{b,t}^{q+} + \Delta D_{b,t}^{q-}\Bigr) \label{eq:cost} \\
    & \hspace{-50pt} \text{subject to} \quad \notag
\end{align}
\textbf{Power Balance at Substation Buses} ($b \in \mathcal{N}^{\text{sub}}$):
\begin{align}
    & p_{b,t}^{\text{sub}} + \sum_{\substack{l \in \mathcal{L} \\ \text{to}(l)=b}} f_{l,t}^{p}
    - \sum_{\substack{l \in \mathcal{L} \\ \text{fr}(l)=b}} f_{l,t}^{p} - D_{b,t}^p \notag \\
    &\hspace{120pt} - \Delta D_{b,t}^{p+} + \Delta D_{b,t}^{p-} = 0\\
    & q_{b,t}^{\text{sub}} + \sum_{\substack{l \in \mathcal{L} \\ \text{to}(l)=b}} f_{l,t}^q
    - \sum_{\substack{l \in \mathcal{L} \\ \text{fr}(l)=b}} f_{l,t}^q
    - D_{b,t}^q \notag \\
    &\hspace{120pt}- \Delta D_{b,t}^{q+} + \Delta D_{b,t}^{q-} = 0
\end{align}

\textbf{Power Balance at Load Buses} ($b \in \mathcal{N} \setminus \mathcal{N}^{\text{sub}}$):
\begin{align}
    & \sum_{\substack{l \in \mathcal{L} \\ \text{to}(l)=b}} f_{l,t}^{p} - \sum_{\substack{l \in \mathcal{L} \\ \text{fr}(l)=b}} f_{l,t}^{p} - D_{b,t}^p - \Delta D_{b,t}^{p+} + \Delta D_{b,t}^{p-} = 0 \\
    & \sum_{\substack{l \in \mathcal{L} \\ \text{to}(l)=b}} f_{l,t}^q
    - \sum_{\substack{l \in \mathcal{L} \\ \text{fr}(l)=b}} f_{l,t}^q
    - D_{b,t}^q - \Delta D_{b,t}^{q+} + \Delta D_{b,t}^{q-} = 0
\end{align}

\textbf{Voltage Equations for Switchable Lines} ($l \in \mathcal{L}^{\text{sw}}$):
\begin{align}
    -v_{b,t}^{\text{fr}(l)} + v_{b,t}^{\text{to}(l)} + 2\!\left(R_l f_{l,t}^{p} + X_l f_{l,t}^q\right) &\leq (1-z_{l,t}^{\text{sw}}) M \\
    v_{b,t}^{\text{fr}(l)} - v_{b,t}^{\text{to}(l)} - 2\!\left(R_l f_{l,t}^{p} + X_l f_{l,t}^q\right) &\leq (1-z_{l,t}^{\text{sw}}) M
\end{align}

\textbf{Voltage Equations for all Lines} ($l \in \mathcal{L}$):
\begin{align}
    -v_{b,t}^{\text{fr}(l)} + v_{b,t}^{\text{to}(l)} + 2\!\left(R_l f_{l,t}^{p} + X_l f_{l,t}^q\right) &\leq (1-av_{l,t}) M \\
    v_{b,t}^{\text{fr}(l)} - v_{b,t}^{\text{to}(l)} - 2\!\left(R_l f_{l,t}^{p} + X_l f_{l,t}^q\right) &\leq (1-av_{l,t}) M
\end{align}

\textbf{Voltage Limits:}
\begin{align}
    v_{b,t} &= V_{\text{ref}}^2 && b \in \mathcal{N}^{\text{sub}} \\
    \underline{V}_b^2 &\leq v_{b,t} \leq \overline{V}_b^2 && b \in \mathcal{N}
\end{align}

\textbf{Thermal Limits for Switchable Lines} ($l \in \mathcal{L}^{\text{sw}}$):
\begin{align}
     -F_l^{\max} \cdot z_{l,t}^{\text{sw}} \leq f_{l,t}^{p} &\leq F_l^{\max} \cdot z_{l,t}^{\text{sw}} \\
    -F_l^{\max} z_{l,t}^{\text{sw}} \leq f_{l,t}^q &\leq F_l^{\max} z_{l,t}^{\text{sw}}
\end{align}

\textbf{Thermal Limits for all Lines} ($l \in \mathcal{L}$):
\begin{align}
    -F_l^{\max} \cdot av_{l,t} \leq f_{l,t}^{p} &\leq F_l^{\max} \cdot av_{l,t}
\end{align}

\textbf{Active and Reactive Power Constraints} ($l \in \mathcal{L}, e \in \{1,2,3,4\}$):
\begin{align}
    f_{l,t}^q &- \cot\!\left[\tfrac{(1/2-e)\pi}{4}\right] (f_{l,t}^{p} - \cos\!\left[\tfrac{e\pi}{4}\right] F_l^{\max}) \leq \sin\!\left[\tfrac{e\pi}{4}\right] F_l^{\max} \\
    -f_{l,t}^q &- \cot\!\left[\tfrac{(1/2-e)\pi}{4}\right] (f_{l,t}^{p} - \cos\!\left[\tfrac{e\pi}{4}\right] F_l^{\max}) \leq \sin\!\left[\tfrac{e\pi}{4}\right] F_l^{\max}
\end{align}

\textbf{Power Injection Limits at Substations:}
\begin{align}
    0 &\leq p_{b,t}^{\text{sub}} \leq P^{\max}_b, && b \in \mathcal{N}^{\text{sub}} \\
    Q^{\min}_b &\leq q_{b,t}^{\text{sub}} \leq Q^{\max}_b && b \in \mathcal{N}^{\text{sub}}
\end{align}

\textbf{Demand Loss and Surplus Bounds:}
\begin{align}
    \Delta D_{b,t}^{p-} &\leq D_{b,t}^p, \quad \Delta D_{b,t}^{q-} \leq D_{b,t}^q && b \in \mathcal{N} \\
    \Delta D_{b,t}^{p\pm}, &\;\Delta D_{b,t}^{q\pm} \geq 0 && b \in \mathcal{N}
\end{align}


\section{Solution Methodology}
\label{sec:methodology}

To solve the problem formulated in the previous section, we propose an approach based on Proximal Policy Optimization (PPO)~\cite{schulman2017proximalpolicyoptimizationalgorithms} with an actor-critic architecture. In the next subsections, we describe the details of our solution methodology.

\subsection{PPO Algorithm}
  
The PPO algorithm tailored for the problem under consideration consists of the following components.

\subsubsection{Actor-Critic Architecture}
 
The framework maintains two separate neural networks: an actor $\pi_\theta: \mathcal{S} \to \mathcal{N}(\boldsymbol{\mu}_\theta(\boldsymbol{s}), \boldsymbol{\sigma}_\theta^2(\boldsymbol{s}))$ that outputs a Gaussian distribution over continuous actions, and a critic $V_\phi: \mathcal{S} \to \R$ that estimates the expected discounted return from a given state. Both networks share a common input representation of the state $\boldsymbol{s}_t$ but have separate parameters $\theta$ and $\phi$. Each network consists of two hidden layers of 256 units with $\tanh$ activations.
 
The actor outputs a mean vector $\boldsymbol{\mu}_\theta(\boldsymbol{s}) \in \R^{d_a}$ and a diagonal covariance $\boldsymbol{\sigma}_\theta^2(\boldsymbol{s}) \in \R^{d_a}_{>0}$, where $d_a$ is the dimension of the output vector. During training, actions are sampled stochastically to encourage exploration; during evaluation, the mean $\boldsymbol{\mu}_\theta(\boldsymbol{s})$ is used deterministically. We project the $d_a$-dimensional sampled action (or the mean vector during evaluation) onto $\{0,1\}^{|\mathcal{L}^\text{sw}|}$ to obtain the switching decisions. Details of the projection scheme are provided in \Cref{sec:action_mapping}.
 
\subsubsection{Training Objective}

PPO optimizes a clipped surrogate objective that constrains how far the updated policy can move from the one used to collect the trajectory. Let
\begin{equation}
    \rho_t(\theta) = \frac{\pi_\theta(\boldsymbol{a}_t \mid \boldsymbol{s}_t)}{\pi_{\theta_{\text{old}}}(\boldsymbol{a}_t \mid \boldsymbol{s}_t)}
\end{equation}
be the probability ratio between the new and old policy. The clipped objective is:
\begin{equation}
    \mathcal{L}^{\text{CLIP}}(\theta) = \E_t\!\left[\min\!\left(\rho_t(\theta)\hat{A}_t,\; \text{clip}\!\left(\rho_t(\theta), 1-\epsilon, 1+\epsilon\right)\hat{A}_t\right)\right]
    \label{eq:ppo_clip}
\end{equation}
where $\epsilon$ is the clip coefficient. By bounding $\rho_t$ within $[1-\epsilon, 1+\epsilon]$, this objective prevents destructively large updates when the advantage estimate is noisy and ensures the new policy stays close to the behavior distribution on which the advantages were computed.

The advantages $\hat{A}_t$ are estimated via Generalized Advantage Estimation (GAE)~\cite{Schulman2015HighDimensionalCC}
The total loss combines the clipped policy objective, an entropy bonus, and a value-function term:
\begin{align}
    \mathcal{L}(\theta, \phi) = &-\left(\mathcal{L}^{\text{CLIP}}(\theta) + c_{\text{ent}} \cdot \mathcal{H}[\pi_\theta(\cdot \mid \boldsymbol{s}_t)]\right) \notag\\
    & + c_{\text{vf}} \cdot \E_t\!\left[\left(V_\phi(\boldsymbol{s}_t) - \hat{G}_t\right)^2\right]
\end{align}
where $\mathcal{H}[\pi_\theta(\cdot\mid\boldsymbol{s}_t)]$ is the policy entropy, $\hat{G}_t = \sum_{t'=t}^{T}\gamma_{RL}^{T-t'}r_{t'}$ is the empirical return, and $c_{\text{ent}}$, $c_{\text{vf}}$ are scaling coefficients. The entropy bonus discourages premature convergence to a deterministic policy, encouraging continued exploration of the topology space during training.
 
Algorithm~\ref{alg:ppo} presents the complete PPO training procedure adapted for the PSPS problem and the key hyperparameters are listed in \Cref{tab:hyperparameters}.
 
\begin{algorithm}[t]
\caption{PPO Training for PSPS with Trajectory Reward Standardization}
\label{alg:ppo}
\begin{algorithmic}[1]
\STATE Initialize actor $\pi_\theta$ and critic $V_\phi$ with random weights
\STATE Pre-solve power flow for all feasible topologies
\FOR{episode $= 1, \ldots, N_{\text{episodes}}$}
    \STATE Reset environment: initialize $\boldsymbol{s}_1$
    \STATE $\mathcal{D} \leftarrow \emptyset$ \hfill \COMMENT{trajectory buffer}
    \FOR{$t = 1, \ldots, T$}
        \STATE Sample action $\boldsymbol{a}_t \sim \pi_\theta(\cdot \mid \boldsymbol{s}_t)$
        \STATE Map $\boldsymbol{a}_t \to \boldsymbol{z}_t^{\text{sw}}$ via discretization (Section~\ref{sec:action_mapping})
        \STATE Compute raw reward $r_t$ via \eqref{eq:reward}
        \STATE Simulate transitions via \eqref{eq:transition} to obtain $\boldsymbol{s}_{t+1}$
        \STATE Store $(\boldsymbol{s}_t, \boldsymbol{a}_t, r_t, \boldsymbol{s}_{t+1}, \log\pi_\theta(\boldsymbol{a}_t|\boldsymbol{s}_t))$ in $\mathcal{D}$
    \ENDFOR
    \STATE Compute trajectory mean and std: $\mu_r \leftarrow \frac{1}{T}\sum_t r_t$,\; $\sigma_r \leftarrow \sqrt{\frac{1}{T}\sum_t (r_t - \mu_r)^2}$
    \STATE Standardize rewards in $\mathcal{D}$: $\tilde{r}_t \leftarrow (r_t - \mu_r)\,/\,(\sigma_r + \varepsilon)$ for all $t$
    \STATE Compute advantages $\hat{A}_t$ using GAE on $\{\tilde{r}_t\}$
    \FOR{update step $= 1, \ldots, K$}
        \STATE Sample mini-batch $\mathcal{B} \subset \mathcal{D}$
        \STATE Update $\theta$, $\phi$ by descending $\nabla_{\theta,\phi} \mathcal{L}(\theta, \phi)$ on $\mathcal{B}$
    \ENDFOR
\ENDFOR
\end{algorithmic}
\end{algorithm}
 
\begin{table}[t]
    \centering
    \caption{PPO Hyperparameters}
    \label{tab:hyperparameters}
    \begin{tabular}{@{}lc@{}}
        \toprule
        Parameter & Value \\
        \midrule
        Learning rate & $4 \times 10^{-4}$ \\
        Entropy coefficient $c_{\text{ent}}$ & 0.05 \\
        Value loss coefficient $c_{\text{vf}}$ & 0.5 \\
        Discount factor $\gamma_{\text{RL}}$ & 0.99 \\
        GAE $\lambda$ & 0.95 \\
        Clip coefficient $\epsilon$ & 0.2 \\
        Training episodes & 10000 \\
        Hidden layer size & 256 \\
        \bottomrule
    \end{tabular}
\end{table}

\subsection{Reward Standardization}
\label{sec:reward_norm}
 
The reward signal \eqref{eq:reward} aggregates the cost of serving the demand of electricity, the penalty for load shedding (or over-supplying) as well as the cost of switching topologies. These three cost components differ substantially in scale and frequency, energy costs are incurred at every step, load loss penalties are large but intermittent, and switching costs are sparse. This produces a heavy-tailed reward distribution within each episode, where time steps involving line failures generate penalties orders of magnitude larger than typical steps. Unnormalized, this heterogeneity couples gradient magnitudes to the reward scale and makes the critic's value targets harder to learn~\cite{vanHasselt2016rewardnormalization}.
 
To address this, we perform per-trajectory reward standardization prior to advantage computation. After collecting the full trajectory $\{r_1, \ldots, r_T\}$, we compute the episode mean and standard deviation:
\begin{equation}
    \mu_r = \frac{1}{T}\sum_{t=1}^{T} r_t, \qquad
    \sigma_r = \sqrt{\frac{1}{T}\sum_{t=1}^{T}(r_t - \mu_r)^2}
\end{equation}
and standardize each reward before passing it to GAE:
\begin{equation}
    \tilde{r}_t = \frac{r_t - \mu_r}{\sigma_r + \varepsilon}
    \label{eq:norm_reward}
\end{equation}
where $\varepsilon$ is a small constant for numerical stability. The resulting $\{\tilde{r}_t\}$ have zero mean and unit variance within every episode, keeping gradient updates and critic targets on a consistent scale regardless of how many failures occurred. A further benefit is hyperparameter portability, the same learning rate and clip coefficient $\epsilon$ transfer across both the 54-bus and 138-bus systems without re-tuning, despite their different absolute cost scales. Reward standardization is applied only during training; raw rewards are used at evaluation time.

\section{Case Study: Effect of Failure Model Nonlinearity on Policy Performance}
\label{sec:case_study}

A central question motivating this work is: How does policy performance change when the true failure model exhibits threshold nonlinearity rather than the linear relationship assumed by existing optimization methods? We investigate this by varying the threshold parameter $\tau$ in the step failure model \eqref{eq:step_model} and measuring performance of all three policies. A low threshold ($\tau = 0.1$) represents a pessimistic scenario where failure risk rises even at moderate flows; a high threshold ($\tau = 0.9$) represents an optimistic one where lines can carry nearly full capacity before risk increases. The key consequence of model mismatch examined here is risk overestimation, when $\tau$ takes intermediate values, the linear model penalizes flows that the true step model considers safe.
\subsection{Experiment Setup}
\begin{figure}[t]
    \centering
    \includegraphics[width=\linewidth]{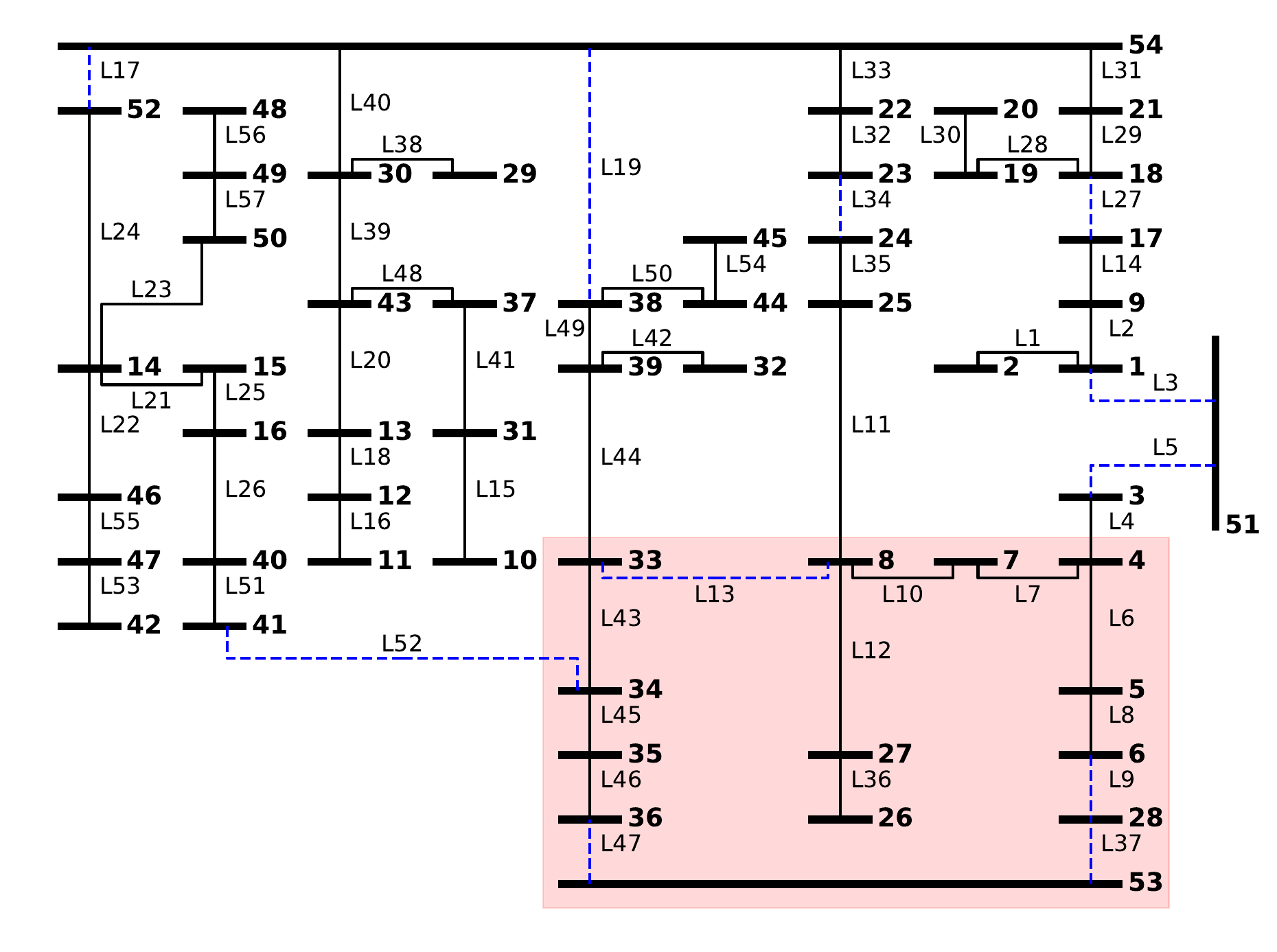}
    \caption{Single-line diagram of the 54-bus distribution system used for sensitivity analysis. The Wildfire-affected and high-threat is highlighted in red. The dashed lines are the switchable power lines.}
    \label{fig:54_bus_system}
\end{figure}

The 54-bus distribution system (\Cref{fig:54_bus_system}) comprises 54 buses (3 substations, 51 load buses), 57 lines (11 switchable), and 14 wildfire-area lines, operated over a 24-hour horizon. The data for the system is taken from \cite{moreira2023dataset}. Cost parameters are: energy price \$10/MWh, load loss cost \$10/MWh, and switching cost \$100/operation. We consider values of thresholds between 0 and 1 at intervals of $0.1$. For each threshold value, we train a dedicated PPO agent for 10000 episodes and evaluate over 200 test episodes under the same true step model. All experiments were conducted on a MacBook Pro with an Apple M4 chip and 16 GB of unified memory.

\subsection{Baseline Policies}
We compare the dynamic PPO-based policy against two static optimization-based policies derived from the outer approximation algorithm of~\cite{moreira2024distribution}. We assume that the baseline risk parameters $\gamma_l$ and $\beta_l$ vary with time in this study with the risk being higher during hours 12 to 20, and for the remaining hours we assume the parameters are only 20\% of the peak values. The fixed baselines below assume the risk parameters are the ones from the peak-hours in order to provide a meaningful benchmark.

\textbf{Static Optimization-Based DIU Policy (Opt-DIU):} This policy is produced by the outer approximation algorithm while ignoring the decision-dependent influence of line flows on failure probabilities (decision-independent uncertainty, DIU). The resulting configuration closes lines $\{3, 5, 37, 47, 52\}$ and opens lines $\{9, 13, 17, 19, 27, 34\}$ throughout the horizon.

\textbf{Static Optimization-Based DDU Policy (Opt-DDU):} This policy accounts for the decision-dependent influence of line flows on failure probabilities but assumes the linear model \eqref{eq:linear_model}. The resulting configuration closes lines $\{3, 17, 19, 34, 37\}$ and opens lines $\{5, 9, 13, 27, 47, 52\}$. This represents the best static policy available from existing linear optimization methods.

Both baselines fix their topology at time $t=0$ and hold it constant throughout the horizon.

\subsection{Handling Topology Constraints and Action Mapping}
\label{sec:action_mapping}

Distribution systems must maintain a radial (tree) structure where all load buses must be connected to exactly one substation path, with no loops. This constrains which combinations of switchable lines can be simultaneously closed.

We identify that the feasible topology space decomposes into five independent constraint groups, as shown in \Cref{tab:topology_groups}. Within each group, at most one line can be closed simultaneously (otherwise a loop would form). The total number of feasible topologies is $3 \times 3 \times 4 \times 3 \times 3 = 324$, compared to $2^{11} = 2048$ without constraints.

\begin{table}[t]
    \centering
    \caption{Decomposition of Topology Constraints (54-bus System)}
    \label{tab:topology_groups}
    \begin{tabular}{@{}clcc@{}}
        \toprule
        Group & Lines & Constraint & Configurations \\
        \midrule
        1 & $\{9, 37\}$ & At most 1 closed & 3 \\
        2 & $\{17, 52\}$ & At most 1 closed & 3 \\
        3 & $\{13, 19, 47\}$ & At most 1 closed & 4 \\
        4 & $\{5, 34\}$ & At most 1 closed & 3 \\
        5 & $\{3, 27\}$ & At most 1 closed & 3 \\
        \midrule
        \multicolumn{3}{l}{Total feasible topologies} & $3 \times 3 \times 4 \times 3 \times 3 = 324$ \\
        \bottomrule
    \end{tabular}
\end{table}

We employ PPO with a continuous action space of dimension 5, where each dimension corresponds to one constraint group. The policy network outputs $\boldsymbol{a} \in \R^5$, which is mapped to discrete topology choices as follows:

\begin{enumerate}
    \item \textbf{Clamp:} $a_i \leftarrow \text{clip}(a_i, -5, 5)$
    \item \textbf{Normalize:} $\tilde{a}_i = \dfrac{a_i + 5}{10} \in [0,1]$
    \item \textbf{Discretize:} $k_i = \min\!\left(\lfloor \tilde{a}_i \cdot n_i \rfloor,\, n_i - 1\right)$
\end{enumerate}
where $n_i$ is the number of feasible configurations for group $i$ and $k_i$ selects which configuration to activate. This procedure guarantees that every output action corresponds to a radially feasible topology, and allows gradient-based optimization via the continuous parameterization.

This procedure generalizes directly to other systems where the constraint groups are identified by determining which sets of switchable lines share a loop, and the output dimension of the actor network equals the number of independent groups. For the 138-bus system considered in Section~\ref{sec:case_study_138}, 12 switchable lines form a different set of constraint groups via the same analysis.

\subsection{Pre-solved Power Flow Solutions}

To accelerate training, we pre-compute optimal power flow solutions for all 324 feasible topologies assuming full line availability. During training episodes with no failures, the appropriate cached solution is retrieved directly. When failures occur, the power flow is re-solved with updated availability constraints. This caching strategy substantially reduces per-episode training time without affecting correctness.

\subsection{Learned Topologies}

The PPO agent's topology selections adapt systematically to $\tau$:

\begin{itemize}
    \item \textbf{Low thresholds ($\tau \in [0.1, 0.3]$):} The agent learns aggressive de-energization strategies, opening more lines in wildfire-prone areas to minimize failure risk. At low thresholds, the step model imposes elevated risk for nearly any positive flow. The PPO agent converges to the same topology as the static DDU baseline, given that this topology is able to aggressively lower the flows in the the wildfire-area without shedding any load.
    \item \textbf{High thresholds ($\tau \in [0.7, 0.9]$):} The agent maintains more lines energized and recognizes that moderate power flows do not significantly increase risk under the step model. In this regime, the PPO policy aligns with the static DIU baseline.
    \item \textbf{Intermediate thresholds ($\tau \in [0.4, 0.6]$):} The PPO agent selects a distinct configuration that differs from both static baselines, reflecting the nuanced trade-off between load-serving capability and failure risk at these threshold levels.
\end{itemize}

\begin{figure*}[t]
    \centering
    \includegraphics[width=\linewidth]{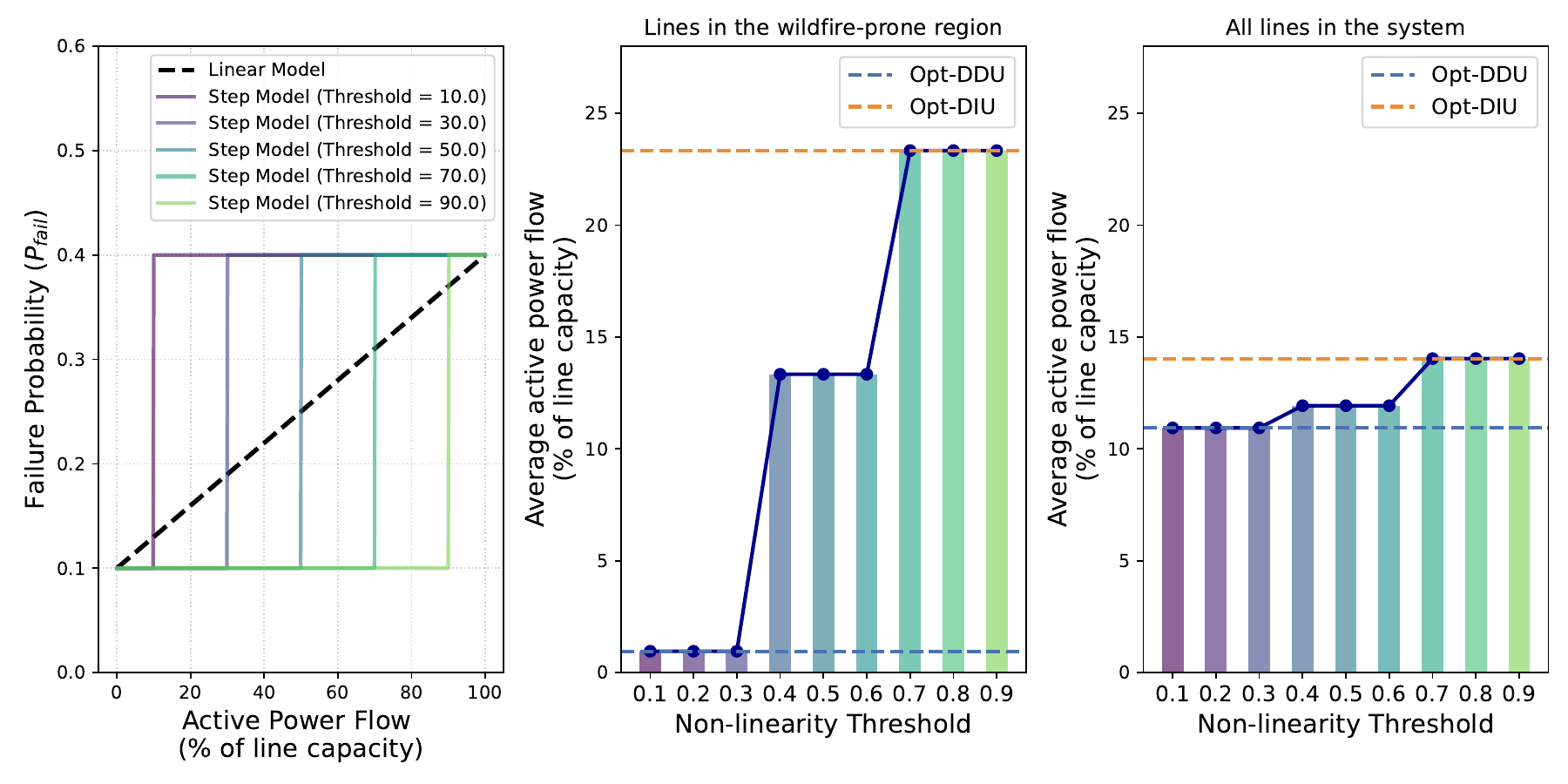}
    \caption{Distribution of power flow magnitudes (as a percentage of line capacity) for wildfire-area lines under the three policies across all threshold values. The PPO-based policy adjusts its loading strategy in response to the observed failure risk.}
    \label{fig:pf_vs_threshold}
\end{figure*}

\subsection{Power Flow Distribution}

\Cref{fig:pf_vs_threshold} shows the power flow distribution for wildfire-area lines across all threshold values. The Opt-DDU policy relies on the linear model \eqref{eq:linear_model}, which assigns nonzero marginal risk to every unit of flow regardless of magnitude. Under the true step threshold model, however, flows below $\tau \cdot F_l^{\max}$ carry no additional failure risk. The linear model therefore overestimates risk in this sub-threshold regime, causing Opt-DDU to de-energize lines and incur switching costs that are simply unnecessary under the true model. The PPO agent, trained directly against the step model, correctly learns that sub-threshold flows are safe and permits higher loading on selected wildfire-area lines as a result. The Opt-DIU policy, which ignores flow-dependent risk entirely, routes substantially more power through the wildfire region, driving its high failure counts

\subsection{Operating Costs and Line Failure Analysis}

\Cref{tab:combined_results} presents average operating costs, switching costs, and line failures per episode for $\tau \in \{0.4, 0.5, 0.6\}$, where the three policies diverge most sharply.

\begin{table}[t]
    \centering
    \caption{Performance by Threshold (54-bus System)}
    \label{tab:combined_results}
    \setlength{\tabcolsep}{4pt}
    \begin{tabular}{@{}llccc@{}}
        \toprule
        & & Dynamic & Static & Static \\
        & & PPO-based & Opt.-based & Opt.-based \\
        $\tau$ & Metric & DDU & DDU & DIU \\
        \midrule
        \multirow{3}{*}{0.4}
            & Op. Cost (\$) & $1696 \pm 612$  & $1912 \pm 57$   & $5869 \pm 1102$ \\
            & Switch Cost (\$) & 400           & 600             & 0               \\
            & Line Failures    & 0.23          & 0.28            & 1.31            \\
        \midrule
        \multirow{3}{*}{0.5}
            & Op. Cost (\$) & $1710 \pm 538$  & $1896 \pm 46$   & $6069 \pm 870$  \\
            & Switch Cost (\$) & 400           & 600             & 0               \\
            & Line Failures    & 0.28          & 0.27            & 1.38            \\
        \midrule
        \multirow{3}{*}{0.6}
            & Op. Cost (\$) & $1703 \pm 624$  & $1842 \pm 61$   & $5333 \pm 1174$ \\
            & Switch Cost (\$) & 400           & 600             & 0               \\
            & Line Failures    & 0.34          & 0.25            & 0.66            \\
        \bottomrule
    \end{tabular}
\end{table}

The PPO agent achieves the lowest total operating cost at all three threshold values. The Opt-DDU policy achieves comparable line failures to PPO at $\tau \in \{0.4, 0.5, 0.6\}$, but does so by over-restricting the network based on the higher risk estimates from the linear model. This overestimation causes Opt-DDU to perform more switching operations than the true failure model warrants (\$600 vs.\ \$400 for PPO), and leads to higher overall operating costs (\$1842--1912 vs.\ \$1696--1710). 

The sharp drop in Opt-DIU line failures from 1.38 at $\tau = 0.5$ to 0.66 at $\tau = 0.6$, a 52\% reduction reveals the step model's structure directly, since the DIU policy routes uniformly high flows regardless of $\tau$, its failure rate tracks the fraction of those flows that exceed the threshold. This also explains why the PPO agent aligns with the Opt-DIU topology for $\tau > 0.6$, where high flows is perceived as safe enough that the cost of aggressive de-energization outweighs its risk benefit.

\section{Case Study: Policy Behavior Under Extreme Wildfire Conditions}
\label{sec:case_study_138}
 
Having established that risk overestimation by the linear model causes Opt-DDU to over-switch under moderate threshold nonlinearity, we now examine the opposite failure mode, risk underestimation. When the true failure model is a step function with a very low threshold and a very high jump in failure probability (something that we can expect during extreme wildfire events), the linear model's gradual slope severely underestimates the risk carried by flows that just exceed the threshold. We use a 138-bus system and fix $\tau = 0.1$ to study this regime, where the consequences of underestimated risk are most pronounced.
 
\subsection{Experiment Setup}
 
\begin{figure*}[t]
    \centering
    \includegraphics[width=0.75\linewidth]{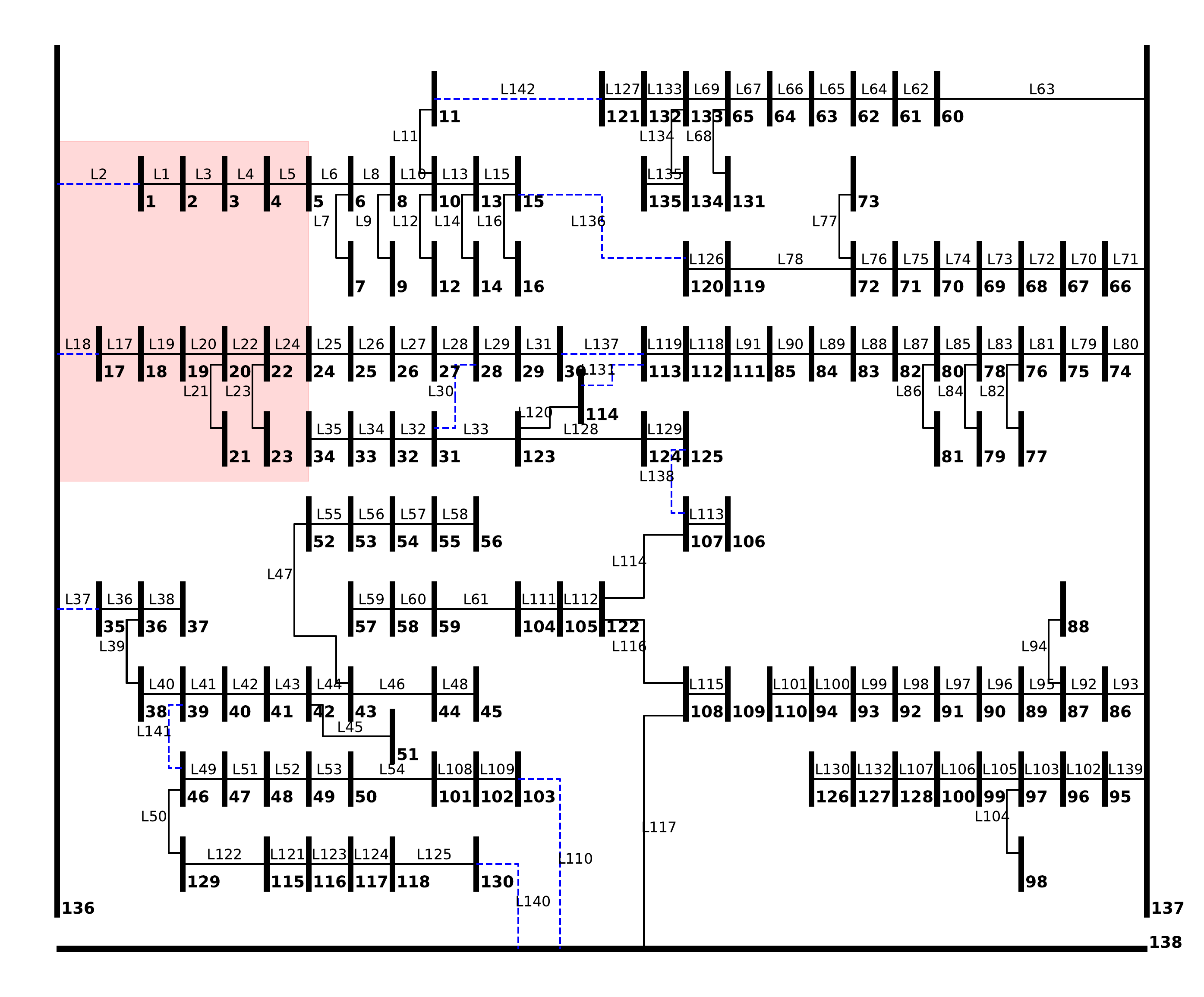}
    \caption{Single-line diagram of the 138-bus distribution system. The Wildfire-affected and high-threat is highlighted in red. The dashed lines are the switchable power lines.}
    \label{fig:138_bus_system}
\end{figure*}

The 138-bus distribution system (\Cref{fig:138_bus_system}) comprises 138 buses (3 substations, 135 load buses), 142 lines (12 switchable), and 13 wildfire-area lines, over a 24-hour horizon. The data for the system is taken from \cite{moreira2023dataset}. Cost parameters are: energy price \$200/MWh, load loss cost \$2000/MWh, switching cost \$200/operation. All experiments were conducted on a MacBook Pro with an Apple M4 chip and 16 GB of unified memory.
 
\textbf{Extreme Wildfire Scenario.} During active wildfire conditions, ambient temperature elevation, reduced wind cooling, and proximity to fire fronts can substantially lower the thermal rating of overhead conductors. We capture this by setting $\tau = 0.1$ and increasing flow sensitivity parameters $\beta_l$ for wildfire-area lines such that the maximum failure probability $\gamma_l + \beta_l \cdot F_l^{\max}$ reaches 90\%. This parameterization closely approximates the conditions utilities face during extreme wildfire events with active fire fronts, where even modest line loading carries significant ignition risk.
 

\subsection{Baseline Policies}
We consider the same two baseline policies considered in the previous case study \Cref{sec:case_study}. For this specific distribution test system, Opt-DIU chooses a topology where the lines \{2, 18, 30, 37, 110\} are closed, while lines \{131, 136, 137, 138, 140, 141, 142\} are kept open. Opt-DDU chooses a topology where the lines \{37, 110, 136, 137, 138\} are closed, while lines \{2, 18, 30, 131, 140, 141, 142\} are kept open. Both these configurations are maintained throughout the horizon.

\subsection{Power Flow Distribution and Load Isolation}
 
\begin{figure}[t]
    \centering
    \includegraphics[width=0.85\linewidth]{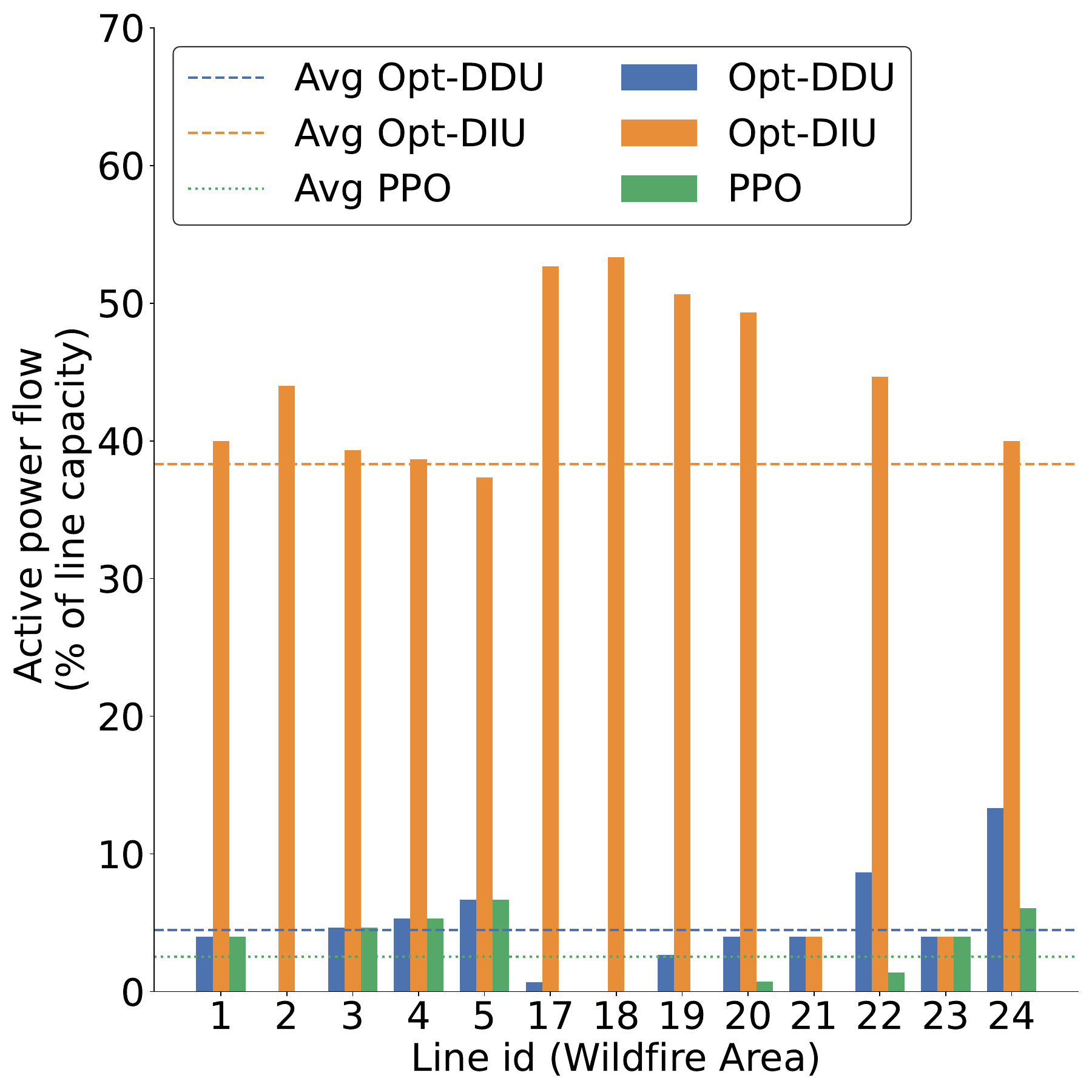}
    \caption{Power flow distribution for wildfire-area lines on the 138-bus system under the extreme scenario ($\tau = 0.1$, max failure probability 90\%). The PPO-based policy drives flows to near zero by isolating loads 17, 18, 19, and 21, while both static baselines continue to route significant power through the high-risk region.}
    \label{fig:power_flows_138}
\end{figure}

\Cref{fig:power_flows_138} compares wildfire-area power flows across the three policies. With $\tau = 0.1$ and maximum failure probability at 90\%, any flow exceeding 10\% of line capacity triggers near-certain failure. The linear model used by Opt-DDU assigns a much more gradual risk increase at these flow levels. Its slope cannot replicate the sharp jump of the step function and therefore underestimates the true risk carried by flows just above the threshold. Consequently, Opt-DDU's optimal topology leaves several wildfire-area lines energized at loadings that appear acceptable under the linear model but are in fact highly dangerous. The PPO agent, trained against the true step model, correctly identifies these flows as unacceptably risky. It responds by de-energizing the supply lines to loads 17, 18, 19, and 21, removing the demand that would otherwise require power to transit the high-risk region and driving wildfire-area flows to near-zero. Similar to what we saw in the previous case study, Opt-DIU routes large amounts of power flow through the wildfire-area because it ignores the flow-dependent risk in it's decision-making.
 
The strategy carries an immediate load shedding penalty at the isolated buses. As \Cref{tab:line_failures_138} shows, however, this deliberate upfront cost is far outweighed by the reduction in uncontrolled, failure-induced load shedding that Opt-DDU and Opt-DIU incur throughout the wildfire region.
 
\subsection{Operating Costs and Line Failures} 
\Cref{tab:line_failures_138} summarizes performance under the extreme scenario. The PPO agent achieves 0.20 failures per episode, compared to 1.14 for Opt-DDU (an 82\% reduction) and 11.06 for Opt-DIU (a 98.2\% reduction). The operating cost advantage is equally significant, PPO's \$327,615 average cost is 17.4\% below Opt-DDU's \$396,996 and 71.2\% below Opt-DIU's \$1,136,818.
 
The Opt-DDU policy performs substantially better than Opt-DIU because its linear-model-derived topology does restrict wildfire-area flows to some degree. However, the linear model's underestimation of risk near $\tau \cdot F_l^{\max}$ means that the restriction is insufficient, flows that Opt-DDU's linear model considers only moderately risky carry a 90\% failure probability under the true step model. The PPO agent learns to eliminate these flows entirely through load isolation rather than merely moderating them.

 Training times as seen in \Cref{tab:training_times} increases by only 34\% from the 54-bus to the 138-bus system, compared to the at least 8X increase in computing times seen in \cite{moreira2024distribution}.

\begin{table}[t]
    \centering
    \caption{Performance Comparison (138-bus System, Extreme Scenario)}
    \label{tab:line_failures_138}
    \resizebox{\columnwidth}{!}{%
    \setlength{\tabcolsep}{4pt}
    \begin{tabular}{@{}lccc@{}}
        \toprule
        & Dynamic & Static & Static \\
        & PPO-based & Opt.-based & Opt.-based \\
        Metric & DDU & DDU & DIU \\
        \midrule
        Op.\ Cost (\$)    & $327{,}615 \pm 42{,}500$ & $396{,}996 \pm 27{,}790$ & $1{,}136{,}818 \pm 15{,}262$ \\
        Switch Cost (\$) & 1200  & 1200  & 0 \\
        Line Failures     & 0.20  & 1.14  & 11.06 \\
        \bottomrule
    \end{tabular}}
\end{table}
 
 
\begin{table}[t]
    \centering
    \caption{Training Time Comparison}
    \label{tab:training_times}
    \begin{tabular}{@{}lcc@{}}
        \toprule
        System & Buses / Lines & Training Time (minutes) \\
        \midrule
        54-bus  & 54 / 57   & 131 \\
        138-bus & 138 / 142 & 176 \\
        \bottomrule
    \end{tabular}
\end{table}
 
\section{Conclusion}
\label{sec:conclusion}
 
This paper presents a reinforcement learning framework for PSPS decision-making that accommodates nonlinear wildfire line failure models. Our PPO-based approach determines switch configurations directly without the need for solving mixed-integer programs. Leveraging per-trajectory reward standardization allows superior training stability, where heterogeneous cost components produce heavy-tailed reward distributions.
 
Sensitivity analysis demonstrates that the PPO agent consistently achieves lower operating costs than both static baselines in the regime where the true step threshold model deviates most from the linear approximation. At intermediate thresholds, the linear model overestimates risk for sub-threshold flows, and therefore the baseline strategies perform switching operations that the true model does not warrant. The PPO agent, trained against the true step model, correctly identifies these flows as safe and avoids this unnecessary cost.
 
Under extreme wildfire threat where the maximum failure probabilities reach close to 1, the linear model severely underestimates the risk of flows at lower levels, and hence, the optimization-based baselines leave wildfire-area lines energized at amounts that are highly dangerous under the true step model. The PPO agent, trained against the true model, correctly identifies these flows as unacceptable and responds with a load isolation strategy that drives wildfire-area flows to near zero. This yields significantly fewer line failures and lower operating costs. 
 
These results establish RL as a promising approach for power system operations under complex, nonlinear risk models. Several directions remain open for future work, including the development of graph neural network policies to enable generalization across network topologies without retraining, the incorporation of probabilistic weather and dynamic line rating forecasts into the state representation to support anticipatory switching, and the replacement of the parametric step threshold model with a failure probability function learned from historical outage records to reduce the need for manual parameter specification.

\section*{Acknowledgment}
This work has been funded by the U.S. Department of Energy, Office of Electricity, under the contract DE-AC02-05CH11231 and National Science Foundation \#2338559. 

\bibliography{scholar}
\bibliographystyle{IEEEtran}

\end{document}